\documentclass[12pt]{article}
\usepackage{latexsym}
\usepackage{amsbsy}
\usepackage{amsmath}   
\usepackage{graphicx}
\usepackage{graphics}
\usepackage{tikz}
\usepackage{pgf}
\usepackage{pgflibraryarrows}
\usepackage{enumerate}
\usepackage{stmaryrd}   

\usepackage{amssymb}

\PassOptionsToPackage{hyphens}{url}
\usepackage{hyperref}

\usetikzlibrary{calc}   

\allowdisplaybreaks[1]   

\newtheorem{thm}{Theorem}

\newtheorem{propn}[thm]{Proposition}

\newcommand{\pf}{\noindent\textit{Proof.\ \ }}
\newcommand{\eopf}{\hfill\hspace*{4pt}\hfill
\fbox{\rule[-1pt]{0pt}{4pt}\hspace*{2pt}}}

\definecolor{greyish}{gray}{0.7}
\definecolor{greenish}{rgb}{0,0.5,0}
\definecolor{lightblue}{rgb}{0.5,0.5,1}
\definecolor{yellowish}{rgb}{1,1,0}
\definecolor{oldgold}{rgb}{0.85,.66,0}
\definecolor{gold}{rgb}{1,0.843,0}
\definecolor{mediumgold}{rgb}{0.75,0.632,0}
\definecolor{darkgold}{rgb}{0.5,0.4215,0}
\definecolor{ochre}{rgb}{0.80,0.47,0.13}
\definecolor{darkbrown}{rgb}{0.5,0.25,0}

\newcommand{\FC}{\hbox{\rm FC}}
\newcommand{\tight}{\hbox{\rm tight}}
\newcommand{\afs}{\hbox{\rm afs}}
\newcommand{\numFC}{\hbox{\rm fcol}}

\newcommand{\dom}{\mathop{\mathrm{dom}}}

\newcommand{\Col}{\mathop{\hbox{Col}}}
\newcommand{\scriptCol}{\mathop{\hbox{\scriptsize Col}}}

\title{The forced colouring function of a graph}
\author{G. E. Farr\thanks{Part of this work was presented at the 31st British Combinatorial Conference, Cardiff, UK, 6--10 July 2026.}   \\
Faculty of I.T.   \\
Monash University   \\
Australia   \\
\href{mailto:Graham.Farr@monash.edu}{\texttt{Graham.Farr@monash.edu}}   \\
\hspace*{1cm}   \\
\tiny
}
\date{15 September 2026}

\begin{document}

\maketitle

\begin{abstract}
The forced colouring function of a graph
gives the probability that a random assignment of colours to a random
subset of vertices can be extended, by a simple local process called \textit{forcing},
to give a proper colouring of the whole graph using the same set of available colours.
This is a polynomial for each fixed number
of colours, and was introduced as a subject for research on the
general theory of graph polynomials.  In this paper we establish its fundamental properties
and give combinatorial interpretations of its derivatives at two particular points.
We also prove
that the problem of computing its value at any specific point in a certain interval is \#P-hard.
\end{abstract}

\section{Introduction}
\label{sec:intro}

Graph colouring is one of the central topics in graph theory, so it is not surprising that
polynomials for counting colourings play a central role in research on graph polynomials.

The first polynomial for counting colourings was the chromatic polynomial $P(G;\lambda)$,
introduced by Birkhoff in \cite{birkhoff1912-1913}.  This gives, for any graph $G$,
the number of $\lambda$-colourings of $G$ for each $\lambda\in\mathbb{N}$.
It helped inspire Whitney \cite{whitney32d} and Tutte \cite{tutte47} to introduce their celebrated bivariate polynomials which bring a surprisingly wide range of
graph-theoretic counting problems into a common framework.  (For the history,
see \cite{farr07a,farr2022,tutte04}.)  Among many other things, the
Tutte-Whitney polynomials can count ``colourings'' according
to their numbers of improperly coloured edges.  Since then, many generalisations of graph colouring
have also led to new graph polynomials for counting them \cite{ashkin-teller1943,averbouch-godlin-makowsky08,averbouch-godlin-makowsky2010,cowen1998,delaharpe-jaeger1995,dohmen-poenitz-tittmann03,harary1985,kotek-makowsky-zilber08,potts52,tittmann2024}.

Graph colouring is also naturally framed as a random process.  Each vertex gets one colour chosen
uniformly at random from a set of $\lambda$ available colours, with the choices for different vertices
being independent.  This does not, in general, give a (proper) $\lambda$-colouring, but we can ask for
the probability that it does.  That probability is given by $P(G;\lambda)/\lambda^n$
(where $n$ is the number of vertices of $G$).

Let us now relax this model by allowing some vertices to remain uncoloured.
The probability of a specific
vertex getting a specific colour is now always $p$,
with $0\le p\le\lambda^{-1}$ and the choices again being independent,
and the probability of a vertex getting no colour at all is $r:=1-\lambda p$.  The outcome will be that
some vertices may be coloured, while others may remain uncoloured,
and there may be edges that are improperly
coloured (in that their endpoints have the same colour).  We can ask for the probability
$\hbox{PC}(G;p,\lambda)$ that the random assignment of colours we get from this process is
proper in the sense that no edge has endpoints that are both coloured with the same colour.
This is the \textit{partial chromatic polynomial} \cite{farr-morgan2025} and is just a very simple algebraic transformation of the generalised chromatic polynomial 
introduced by Dohmen, P\"onitz and Tittmann \cite{dohmen-poenitz-tittmann03}.

Many graph colouring algorithms are based on extending partial colourings to the entire graph while
not using too many extra colours \cite{biggs1990,husfeldt2015,welsh-powell1967}.
So, given a partial colouring, it is of interest to know if it can be extended to a colouring of the entire
graph while using the same set of available colours.
The probability that this is possible, under our random partial colouring model, is denoted by
$\hbox{EC}(G;p,\lambda)$ and called the \textit{extendable colouring function} \cite{farr-morgan2025}.

We can ask for more.  
What is the probability that the random partial assignment of colours can be extended \textit{uniquely}
to the entire graph (again using the same colour set)?  This is $\hbox{UC}(G;p,\lambda)$,
the \textit{uniquely extendable colouring function}.

One way in which a partial assignment of colours can have a unique extension that colours
the whole graph is if the extension is \textit{forced} in the following sense.
Whenever we have an uncoloured vertex $v$ for which
all colours except one appear among its neighbours,
$v$ is given that one remaining available colour.
This is because, in such circumstances, any extension to the whole graph must assign that
particular colour to $v$.  We keep doing this forcing for as long as possible.
(This was studied for $\lambda=3$ in \cite{farr94}.)
We write $\hbox{FC}(G;p,\lambda)$ for the probability that a random partial assignment of
$\lambda$ colours forces, in this sense, a $\lambda$-colouring of $G$ \cite{farr-morgan2025}.
If a $\lambda$-colouring of $G$ is indeed forced, then the partial assignment we started with is
uniquely extendable to a $\lambda$-colouring of $G$, although the converse is not true in general.

The functions PC, EC and FC were introduced in \cite{farr-morgan2025}.  The motivation there was
to provide some intriguing ``specimens'' for study by researchers who are working towards
a ``zoology'' of graph polynomials.  That research programme, led by Makowsky and colleagues
\cite{ellis-monaghan-etal2016,ellis-monaghan-etal2019,godlin-katz-makowsky2012,kotek-makowsky-ravve2013,kotek-makowsky-zilber08,makowsky06,makowsky08,makowsky2012,makowsky2023,makowsky-ravve2016,makowsky-ravve-kotek2019}, has shed much light on graph polynomials as well as raising interesting questions
about them.

The functions EC, UC, FC are in some ways unlike other polynomials considered in that research.
This is partly because they involve counting structures that satisfy
conditions that are more complex than those involved in many classical graph polynomials.
For example, colourings, flows and stable sets are determined just by conjunctions of simple
local conditions based on vertices, edges or neighbourhoods, but the conditions that must be
satisfied by a partial colouring in order for it to be extendable in these various ways are significantly
more complex \cite{farr-morgan2025}.

These new polynomials raise, and illustrate, some issues that are worth consideration in the development of a
general theory.  For example, many graph polynomials can be expressed as a sum over subsets
of the edge set (or vertex set) of the graph, and many also satisfy \textit{reduction relations},
i.e., simple linear recurrence relations involving a fixed finite number of local graph operations
(such as edge deletion, edge contraction, vertex deletion, neighbourhood deletion, etc.).
In fact, it has been shown by Godlin, Katz and Makowsky \cite{godlin-katz-makowsky2012}
that all graph polynomials from a very general class that satisfy a reduction relation must also be
expressible as a sum over subsets.  But not every graph polynomial has a reduction relation,
so we cannot necessarily go the other way.  Nonetheless, it turns out that graph polynomials that
do not satisfy a reduction relation in their original class (usually, graphs) often \textit{do} satisfy
such a relation in some more general class of combinatorial objects.  This is discussed at length
in \cite[\S5]{farr-morgan2025}, where many examples are given, and the forced colouring
polynomial $\FC$, treated for each fixed $\lambda\in\mathbb{N}$ as a univariate
polynomial in $p$, is such a case.  We do not yet know how widespread this intriguing
phenomenon is.  Could a kind of converse of the Godlin-Katz-Makowsky result be proved?
This would entail showing that
every graph polynomial (of some general type) that admits a sum-over-subsets expression
has a reduction relation over some class
of combinatorial objects that may be wider than just graphs.
This question and some related ones are raised and discussed in \cite[\S\S5--6]{farr-morgan2025}.
(See also \cite[p.~574]{tittmann2025}.)
We do not address those general questions here, but they motivate our interest in the forced colouring
function as a kind of test case for some of these questions.
So it is to be hoped that better understanding
of the forced colouring function and its relatives may contribute to some improved understanding
of graph polynomials in general.

This paper focuses on the forced colouring function FC, since in some respects it is more
tractable than EC or UC .  After giving definitions and notation in \S\ref{sec:defn-notn},
we establish basic properties of $\FC$ in \S\ref{sec:basic-properties}.  These include the interpretation
of its value at $p=(\lambda+1)^{-1}$,
some special forms the polynomial has for certain values of $\lambda$, a link with
the chromatic polynomial, some simple inequalities, an expression for bipartite graphs when
$\lambda=2$, some examples, and some basic algebraic observations.
Then, in \S\ref{sec:derivatives}, we give combinatorial interpretations of the derivatives at
$p=(\lambda+1)^{-1}$ and $p=\lambda^{-1}$.  Finally, in \S\ref{sec:complexity}, we consider its
computational complexity.  We prove that
finding its value for any fixed $\lambda\ge3$ and any fixed $p\in(0,\lambda^{-1}]$ is \#P-hard,
which is evidence of intractability.
We close in \S\ref{sec:future-work} with some suggestions for future research.

\section{Definitions and notation}
\label{sec:defn-notn}

Throughout, $G=(V,E)$ is a graph with $n$ vertices and $m$ edges.  The set of vertices that
are adjacent to $v\in V$ in $G$ is denoted by $N_G(v)$.  The maximum degree of a vertex in $G$
is denoted by $\Delta(G)$.
The number of components of $G$ is denoted by $k(G)$.
If $X\subseteq E$ then $V(X)$ denotes the set of vertices of $G$ that are incident with at least
one edge in $X$ (overloading $V$ slightly).
If $U\subseteq V$ then $G[U]$ is the subgraph of $G$ induced by $U$.

We write $[k] = \{1,2, \ldots, k\}$.

The falling factorial $\lambda(\lambda-1)\cdots(\lambda-k+1)$ is denoted by $(\lambda)_k$.

The restriction of a function $g$ to a subset $Y$ of its domain is denoted by $\left.g\right|_Y$.
Rest assured that this is seldom used in contexts where set cardinality notation is also needed!

We sometimes use the Iverson bracket.  If $Q$ is a proposition then
\[
\llbracket\hbox{$Q$}\,\rrbracket ~=~
\left\{
\begin{array}{cl}
1,  &  \hbox{if $Q$ is true};   \\
0,  &  \hbox{if $Q$ is false}.
\end{array}
\right.
\]

\subsubsection*{Partial assignments, partial colourings and extensions}

Let $\Lambda$ be a finite set.  We think of it as the set of available colours,
and mostly we use $\Lambda=[\lambda]$ where $\lambda\in\mathbb{N}$.

A \textit{partial $\Lambda$-assignment} is a function $f:W\rightarrow\Lambda$ where
$W\subseteq V$.  If $\lambda\in\mathbb{N}$ and $\Lambda=[\lambda]$ then we call
$f$ a \textit{partial $\lambda$-assignment}.  If such an assignment has domain $V$ then
we may call it a \textit{total $\Lambda$-assignment} or just a \textit{$\Lambda$-assignment}
(and similarly with $\lambda$ replacing $\Lambda$ when $\lambda\in\mathbb{N}$).
The set of partial $\Lambda$-assignments with domain $W$ is just the set $\Lambda^W$,
and we use this compact notation sometimes.

For every $i\in[\lambda]$, its \textit{colour class} $C(i)=C_f(i)$ under
a partial $\lambda$-assignment $f$ of $G$ is given by
$C(i)=f^{-1}(i)=\{v\in V:f(v)=i\}$.
Every colour class $C(i)$ induces a subgraph $G[C(i)]$ of~$G$.

A colour class is \textit{stable} if no two of its vertices are adjacent, i.e., it is a stable set (or independent set) in $G$, otherwise it is \textit{unstable}.
A partial $\lambda$-assignment $f$ is \textit{stable} if
every colour class is stable, and it is then a \textit{partial $\lambda$-colouring}.
Note that a partial $\lambda$-assignment $f$ is stable if and only if it is
a colouring of $G[\dom f]$.\footnote{We mostly avoid referring to ``proper'' colourings from now on, since (a) colourings are proper by definition, and (b) we use ``proper'' terminology for extensions of partial $\lambda$-assignments in a way that may create confusion if we use that term for partial $\lambda$-assignments/colourings too.}
If a partial $\lambda$-assignment is not a partial $\lambda$-colouring then
it is \textit{unstable}, since then it has at least one unstable colour class.
A partial $\lambda$-colouring with domain $V$ is just a $\lambda$-colouring.
For any graph $G$ and positive integer $\lambda$, the set of all $\lambda$-colourings of $G$ is denoted by
$\Col(G;\lambda)$.

An \textit{extension} of a partial $\lambda$-assignment $f$
is a partial $\lambda$-assignment $g$ such that $\dom f\subseteq\dom g$ and
$f(v)=g(v)$ for all $v\in\dom f$.
The extension is \textit{proper} if $\dom f$ is a proper subset of $\dom g$.
(These terms are just instances of standard function terminology in this context.)
This terminology does not require either $f$ or $g$ to be a partial $\lambda$-colouring.

If some proper extension of $f$ is actually a partial $\lambda$-colouring of $G$
then we say $f$ is \textit{properly $\lambda$-extendable},
or just \textit{properly extendable} if $\lambda$ is clear from the context.
For this to happen, it is necessary (but not sufficient) that $f$ be a partial $\lambda$-colouring of~$G$.
In general, a partial $\lambda$-colouring $f$ may or may not be extendable
to a $\lambda$-colouring of~$G$.
If $f$ is not extendable to a $\lambda$-colouring of $G$ then we say it is $\lambda$\textit{-contradictory},
or just \textit{contradictory} if $\lambda$ is clear from the context.

\subsubsection*{Forcing}

Let $f$ be a partial $\lambda$-assignment.
A vertex $v\not\in\dom f$ is \textit{immediately $\lambda$-forced by $f$} if the number of
different colours
appearing among its neighbours is $\lambda-1$, i.e.,
\[
\left| \{f(v)\mid v\in N_G(v)\cap\dom f\} \right| = \lambda-1 .
\]
In this case, every extension of $f$ that is a partial $\lambda$-colouring and includes $v$ in its
domain must give $v$ the sole colour that does not appear among its neighbours.  
We write $f;v$ for the unique extension of $f$ to $(\dom f)\cup\{v\}$ that does this:
\begin{eqnarray*}
f;v \,:\, (\dom f)\cup\{v\} & \rightarrow & [\lambda] ,   \\
\left.(f;v)\right|_{\dom f} & = & f ,   \\
(f;v)(v) & = & \hbox{the sole member of $[\lambda]\setminus\{f(v)\mid v\in N_G(v)\cap\dom f\}$}.
\end{eqnarray*}
If the number of different colours that appear among the neighbours of $v$ is $\lambda$, then
none of our $\lambda$ colours can be given to $v$ without creating bad edges.
On the other hand, if that number is $\le\lambda-2$,
then the colour of $v$ is not forced by the colours of its neighbours, as there are at least
two options for it.

The vertex $v$ is \textit{(eventually) $\lambda$-forced by $f$} if there is a sequence
of distinct vertices $v_1,\ldots,v_k=v\in V\setminus\dom f$,
and a sequence of partial $\lambda$-assignments $f=f_0,f_1,\ldots,f_k$ such that,
for all $i\in\{1,\ldots,k\}$,
\begin{itemize}
\item $v_i$ is immediately forced by $f_{i-1}$,
\item $f_i=f_{i-1};v_i$.
\end{itemize}
We also say that each $f_i$ is \textit{$\lambda$-forced by $f$}, where $0\le i\le k$.

We write $\Phi^*f$ for an extension of $f$ that is constructed from $f$ by continuing
to force vertices for as long as possible, assigning to each forced vertex the sole colour that it
is forced to have.  Formally, $\Phi^*f$ is defined by the following conditions:
\begin{itemize}
\item there exists a sequence of distinct vertices $v_1,\ldots,v_k\in V\setminus\dom f$ such that
$\Phi^*f=f;v_1;v_2;\cdots;v_k$, and
\item there is no vertex $w\in V\setminus(\dom f\cup\{v_1,\ldots,v_k\})$ such that $\Phi^*f$ forces~$w$.
\end{itemize}
If $f$ is extendable to a $\lambda$-colouring of $G$ then $\Phi^*f$ is unique, although it might not be
a $\lambda$-colouring of $G$.  If $f$ is contradictory then $\Phi^*f$ may or may not be unique.

Throughout, we may drop $\lambda$ from ``$\lambda$-forced'' when it is clear from the context.

Note that if $f$ is unstable then it is not properly $\lambda$-extendable
and it cannot force a $\lambda$-colouring of $G$.

The following basic observation should also be kept in mind.

\begin{propn}
If $f$ is not total and
$\lambda$-forces a $\lambda$-colouring of $G$ and $\mu>\lambda$ then $f$ does
not $\mu$-force a $\mu$-colouring of $G$.
\end{propn}

\pf
Let $v$ be any vertex not in $\dom f$.  (Such a vertex exists, since $f$ is not total.)
The number of different colours appearing among the neighbours of $v$ cannot be $\lambda$,
because then $v$ could not be coloured, whereas we know $f$ forces a $\lambda$-colouring
of $G$.  So this number of neighbouring colours is $\le\lambda-1<\mu-1$, so the neighbours of $v$
do not have enough colours to immediately $\mu$-force $v$.  Therefore no uncoloured vertex
is immediately $\mu$-forced by $f$.  Therefore $f$ does not $\mu$-force any other vertices
and (since $f$ is not total) it does not $\mu$-force a $\mu$-colouring of $G$.
\eopf   \\

This does not preclude the possibility that some extension of $f$ might $(\lambda+1)$-force
a $(\lambda+1)$-colouring of $G$.

Forcing was considered in \cite{farr94}, where it was shown that determining, for a given input graph,
whether or not a given partial 3-colouring forces a 3-colouring of the entire graph is logspace-complete for P.
In the present work, we are interested in what gets forced by random partial colourings.

\subsubsection*{Random colouring and associated functions}

Given a graph $G$, a positive integer $\lambda$ representing some number of available
colours, and a probability $p\le\lambda^{-1}$, put $r:=1-\lambda p$ and define the following
random colouring model, noting that the ``colouring'' it produces may be improper (i.e., unstable).

Each vertex $v\in V(G)$ remains uncoloured with probability $r$ and otherwise is given a colour
chosen uniformly at random from the $\lambda$ available colours.  So, for any specific colour,
the probability that $v$ gets that colour is $p$.  The choices made at different vertices are
independent.  This process generates a random partial $\lambda$-assignment whose
domain is a subset of $V(G)$.

For every graph $G$ we define four bivariate functions based on partial colourings.

The \textit{partial chromatic polynomial} $\hbox{PC}(G;p,\lambda)$ is defined by
\[
\hbox{PC}(G;p,\lambda)  =  \Pr(\hbox{$f$ is a $\lambda$-colouring of $G[\dom f]$}).
\]
This is a simple algebraic transformation of
the \textit{generalised chromatic polynomial} $P(G;x,y)$ introduced
by~\cite{dohmen-poenitz-tittmann03}.  We discuss it and its close relatives briefly
in \cite{farr-morgan2025}.

The \textit{extendable colouring function} is given by
\[
\hbox{EC}(G;p,\lambda)  =  \Pr(\hbox{$f$ has an extension that is a $\lambda$-colouring of $G$}) .
\]
This is discussed in \cite{farr-morgan2025} where it is proposed as another subject of study
in the development of the theory of graph polynomials.

The \textit{uniquely extendable colouring function} is given by
\[
\hbox{UC}(G;p,\lambda)  =  \Pr(\hbox{$f$ has a \textit{unique} extension that is a $\lambda$-colouring of $G$}) .
\]

The \textit{chromatic forcing function} is given by
\[
\hbox{FC}(G;p,\lambda)  =  \Pr(\hbox{$f$ eventually forces a $\lambda$-colouring of $G$}) .
\]
This was introduced in \cite{farr-morgan2025} and is the topic of this paper.

\section{Basic properties}
\label{sec:basic-properties}

Define $\numFC(G,C;\lambda)$ to be the number of partial $\lambda$-assignments $f$ such that
$\dom f=C$ and $f$ eventually forces a $\lambda$-colouring of $G$.  Any partial $\lambda$-assignment
counted by $\numFC(G,C;\lambda)$ will in fact be a partial $\lambda$-colouring.
Observe that
\begin{eqnarray}
\label{eq:numFC-at-empty}
\numFC(G,\emptyset;\lambda)  & = & 
\left\{
\begin{array}{cl}
2^n,  &  \hbox{if $\lambda=1$ and $G$ has no edges};   \\
0,  &  \hbox{otherwise}.
\end{array}
\right.   \\
\label{eq:numFC-at-V}
\numFC(G,V;\lambda)  & = &  P(G;\lambda) .
\end{eqnarray}
The first of these observations
is because an uncoloured isolated vertex is forced when $\lambda=1$ but not otherwise,
and the second is because a partial $\lambda$-assignment with domain $V$ that forces a
$\lambda$-colouring must itself be a $\lambda$-colouring.

By partitioning the sample space according to the domain of the random partial $\lambda$-assignment,
we have

\begin{propn}{\cite[p.~279]{farr-morgan2025}}
\label{propn:FC-sum-over-C}
\begin{equation}
\label{eq:FC-sum-over-C}
\FC(G;p,\lambda)  =  \sum_{C\subseteq V} \numFC(G,C;\lambda) p^{|C|}(1-\lambda p)^{n-|C|} .
\end{equation}
\eopf
\end{propn}
This shows that, for each fixed $\lambda\in\mathbb{N}$, the forced colouring function $\FC(G;p,\lambda)$
is actually a polynomial in $p$.  But it is not, in general, a polynomial in $\lambda$.
Our focus will be on what happens for fixed $\lambda\in\mathbb{N}$.  To emphasise this
perspective we will write
\[
\FC_{\lambda}(G;p) ~:=~ \FC(G;p,\lambda)
\]
and when $\lambda$ is fixed we will refer to $\FC_{\lambda}(G;p)$, a univariate polynomial in $p$,
as the \textit{forced $\lambda$-colouring polynomial} of $G$, dropping $\lambda$ when it is clear
from the context.

We had better note, in passing, the natural connection between $\FC_{\lambda}(G;p)$ and its
related generating polynomial.

For each $i\in\mathbb{N}\cup\{0\}$,
define $\numFC(G,i;\lambda)$ to be the number of partial $\lambda$-assignments $f$ of $G$ such that
$|\dom f|=i$ and $f$ forces a $\lambda$-colouring of $G$.  So
\[
\numFC(G,i;\lambda) ~=~ \sum_{C\subseteq V:|C|=i} \numFC(G,C;\lambda) .
\]
Write $\hbox{FCGP}_{\lambda}(G,x)$ for the generating polynomial for the $\hbox{FC}(G,i;\lambda)$:
\[
\hbox{FCGP}_{\lambda}(G,x)
= \sum_{C\subseteq V} \numFC(G,C;\lambda) x^{|C|}
= \sum_{i=0}^n \numFC(G,i;\lambda) x^i .
\]

\begin{propn}
\begin{eqnarray*}
\FC_{\lambda}(G;p)
& = &
(1-\lambda p)^{n} \, \hbox{\rm FCGP}_{\lambda}\Bigl( G, \frac{p}{1-\lambda p} \,\Bigr) ,   \\
\hbox{\rm FCGP}_{\lambda}(G,x)  & = &
(1+\lambda x)^{n} \, \FC_{\lambda}\Bigl( G; \frac{x}{1+\lambda x} \,\Bigr) .
\end{eqnarray*}
\eopf
\end{propn}

We briefly record some elementary properties of $\FC_{\lambda}(G;p)$.

\begin{propn}
\label{propn:elem-properties}

(a)  If $\lambda<\chi(G)$ then $\FC_{\lambda}(G;p)\equiv0$.

(b)  If $G$ has at least one edge or $\lambda\ge2$ then $\FC_{\lambda}(G;0)=0$.

(c)  $\lambda^{n}\FC_{\lambda}(G;\lambda^{-1}) = P(G;\lambda)$.

(d)  $(\lambda+1)^n\FC_{\lambda}(G;(\lambda+1)^{-1})$ is the number of partial $\lambda$-assignments of $G$ that force a $\lambda$-colouring of $G$.

(e)  For each $\lambda\in\mathbb{N}$, the polynomial $\FC_{\lambda}$ is multiplicative over disjoint
unions:
\[
\FC_{\lambda}(G\cup H;p) ~=~ \FC_{\lambda}(G;p) \, \FC_{\lambda}(H;p) .
\]
\end{propn}

\pf
(a)  If $\lambda<\chi(G)$ then $G$ has no $\lambda$-colouring, so no partial $\lambda$-assignment
can be extended to a $\lambda$-colouring.

(b)  If $p=0$ then, in this case, the random partial $\lambda$-assignment has
domain $\emptyset$, and $\numFC(G,\emptyset;\lambda)=0$ by \ref{eq:numFC-at-empty}.

(c)
If $p=\lambda^{-1}$ then there are no uncoloured vertices, so the random partial $\lambda$-assignment
forces a $\lambda$-colouring of $G$ if and only if it is already a colouring of $G$.

(d)
When $p=(\lambda+1)^{-1}$,
\begin{equation}
\label{eq:prob-factor-with-one-on-lambda-plus1}
p^{|C|}(1-\lambda p)^{n-|C|} = \left(\frac{1}{\lambda+1}\right)^n .
\end{equation}
Using (\ref{eq:FC-sum-over-C}) and (\ref{eq:prob-factor-with-one-on-lambda-plus1}),
\[
(\lambda+1)^n\FC_{\lambda}(G;(\lambda+1)^{-1})  =  (\lambda+1)^n \sum_{C\subseteq V} \numFC(G,C;\lambda) \left(\frac{1}{\lambda+1}\right)^n
= \sum_{C\subseteq V} \numFC(G,C;\lambda).
\]

(e)  Routine.
\eopf   \\

We now establish a simple link with the chromatic polynomial.

\begin{propn}
\label{propn:FC-and-chrompoly}
For every graph $G$ and every $\lambda\in\mathbb{N}$,
\begin{equation}
\label{eq:FC-ge-chrompoly-times-p-to-n}
\FC_{\lambda}(G;p) \ge P(G;\lambda)\,p^n .
\end{equation}
Furthermore, if $\lambda>\Delta(G)+1$ then
\begin{equation}
\label{eq:FC-eq-chrompoly-times-p-to-n}
\FC_{\lambda}(G;p) = P(G;\lambda)\,p^n .
\end{equation}
\end{propn}

\pf
\begin{eqnarray*}
\FC(G;p,\lambda)  & = &
\sum_{C\subseteq V} \numFC(G,C;\lambda) \, p^{|C|}(1-\lambda p)^{n-|C|}
~~~~ \hbox{(by (\ref{eq:FC-sum-over-C}))}   \\
& \ge &
\numFC(G,V;\lambda) p^n
~~~~ \hbox{(just taking the term for $C=V$)}   \\
& = &
P(G;\lambda)\, p^n
~~~~ \hbox{(by (\ref{eq:numFC-at-V}))}.
\end{eqnarray*}

Suppose $\lambda>\Delta(G)+1$.

Let $v$ be any vertex of $G$.  Since $\deg(v)\le\Delta(G)<\lambda-1$, the number of different colours
on the neighbours of $v$ is always $<\lambda-1$, which is never enough to force the colour of $v$.
So the only way to force $v$ is to colour it initially rather than wait for the forcing process to reach it.
Since this holds for every vertex, the initial random partial $\lambda$-assignment must in fact colour every
vertex of $G$, and then the only way for it to \textit{force} a $\lambda$-colouring is for it to
\textit{be} a $\lambda$-colouring.
For every $\lambda$-colouring of $G$, the probability that this is the initial partial $\lambda$-assignment chosen by our random process is $p^n$, so $\FC(G;p,\lambda)=P(G;\lambda)p^n$.
\eopf   \\

Let us summarise the inequalities we have noted so far.
\begin{equation}
\label{eq:ineqs-for-polys}
P(G;\lambda)p^n \le \hbox{FC}(G;p,\lambda) \le \hbox{UC}(G;p,\lambda) \le \hbox{EC}(G;p,\lambda) \le \hbox{PC}(G;p,\lambda)
\end{equation}
whenever $\lambda\in\mathbb{N}$ and $0\le p\le \lambda^{-1}$.  The first is
just (\ref{eq:FC-ge-chrompoly-times-p-to-n}) and the others follow from the logical implications among
the conditions the partial $\lambda$-colouring must satisfy, as discussed in \S\ref{sec:intro}
and \cite[p.~279]{farr-morgan2025}.

Given that the forced colouring function $\FC(G;p,\lambda)$ of a graph $G$
is not a bivariate polynomial in $p$ and $\lambda$, it is not
immediately obvious how to represent it in general.
The numbers $\numFC(G,k;\lambda)$ form an infinite set,
since they depend on $\lambda\in\mathbb{N}$ as well as on $k\in[n]$,
so (\ref{eq:FC-sum-over-C}) in Proposition \ref{propn:FC-sum-over-C} does not immediately
yield a finite representation of $\FC(G;p,\lambda)$.
But (\ref{eq:FC-eq-chrompoly-times-p-to-n}) in Proposition \ref{propn:FC-and-chrompoly}
tells us that, if $\lambda>\Delta(G)+1$, we really just need to state the chromatic polynomial
(along with the simple factor $p^n$),
and stating this once gives a compact description of $\FC(G;p,\lambda)$ for all $\lambda>\Delta(G)+1$.
So it remains to represent $\FC(G;p,\lambda)$ for $\lambda\in[\Delta(G)+1]$.
This can be done, in principle, by just stating all the numbers $\numFC(G,k;\lambda)$ for
$k\in[n]$ and $\lambda\in[\Delta(G)+1]$, and there are $n(\Delta(G)+1)\le n^2$ of these numbers.
So this gives a finite, and reasonably compact, representation of $\FC(G;p,\lambda)$.

The question of how to actually \textit{compute} these functions is a different one.
We consider this in \S\ref{sec:complexity}.  For now, we just give the results of some very
simple computations for a few small examples.

We can use Propositions \ref{propn:FC-sum-over-C}, \ref{propn:elem-properties}(a) and \ref{propn:FC-and-chrompoly}
to compute $\FC_{\lambda}(G;p)$ for some simple graphs.
We do this now for empty graphs, complete graphs, and
all connected graphs of at most four vertices.  (Disconnected graphs
can be handled with Proposition \ref{propn:elem-properties}(e).)  We omit the routine proofs.

\newcommand{\pathThree}{
\begin{tikzpicture}
[vertex/.style={circle,fill,inner sep=0pt,minimum size=0.1cm}]
\node[vertex,fill=violet,draw=black,thick] (a) at (0,0) {};
\node[vertex,fill=violet,draw=black,thick] (b) at (0.3,0) {};
\node[vertex,fill=violet,draw=black,thick] (c) at (0.6,0) {};
\node[vertex,fill=violet,draw=black,thick] (d) at (0.9,0) {};
\draw (a) -- (b);
\draw (b) -- (c);
\draw (c) -- (d);
\end{tikzpicture}
}
\newcommand{\threePan}{
\begin{tikzpicture}
[vertex/.style={circle,fill,inner sep=0pt,minimum size=0.1cm}]
\node[vertex,fill=violet,draw=black,thick] (a) at (0,0.15) {};
\node[vertex,fill=violet,draw=black,thick] (b) at (0,-0.15) {};
\node[vertex,fill=violet,draw=black,thick] (c) at (0.2598,0) {};   
\node[vertex,fill=violet,draw=black,thick] (d) at (0.5598,0) {};
\draw (a) -- (b);
\draw (b) -- (c);
\draw (a) -- (c);
\draw (c) -- (d);
\end{tikzpicture}
}
\newcommand{\WelshGraph}{
\begin{tikzpicture}
[vertex/.style={circle,fill,inner sep=0pt,minimum size=0.1cm}]
\node[vertex,fill=violet,draw=black,thick] (a) at (0,0.15) {};
\node[vertex,fill=violet,draw=black,thick] (b) at (0,-0.15) {};
\node[vertex,fill=violet,draw=black,thick] (c) at (0.3,0.15) {};
\node[vertex,fill=violet,draw=black,thick] (d) at (0.3,-0.15) {};
\draw (a) -- (b);
\draw (a) -- (c);
\draw (c) -- (d);
\draw (b) -- (d);
\draw (b) -- (c);
\end{tikzpicture}
}

\begin{propn}
\label{propn:FC-some-simple-gphs}
For each $\lambda\in\mathbb{N}$ and $n\in\mathbb{N}$,   \\
\begin{eqnarray}
\FC_{\lambda}(\overline{K_n};p)
& = &
\left\{
\begin{array}{ll}
1,  &  \hbox{if $\lambda=1$};   \\
(\lambda p)^n,  &  \hbox{if $\lambda\ge2$}.
\end{array}
\right.   \\[5pt]
\FC_{\lambda}(K_n;p)
& = &
\left\{
\begin{array}{ll}
0,  &  \hbox{if $\lambda<n$};   \\
n!p^{n-1}(n-(n^2-1)p),  &  \hbox{if $\lambda=n$};   \\
\lambda_{(n)} p^n,  &  \hbox{if $\lambda>n$}.
\end{array}
\right.   \\[5pt]
\FC_{\lambda}(K_{1,2};p)
& = &
\left\{
\begin{array}{ll}
0,  &  \hbox{if $\lambda=1$};   \\
2p(3-9p+7p^2),  &  \hbox{if $\lambda=2$};   \\
6p^2(1-2p),  &  \hbox{if $\lambda=3$};   \\
\lambda(\lambda-1)^2p^3,  &  \hbox{if $\lambda\ge4$}.
\end{array}
\right.   \label{eq:FC-K12}   \\[5pt]
\FC_{\lambda}(K_{1,3};p)
& = &
\left\{
\begin{array}{ll}
0,  &  \hbox{if $\lambda=1$};   \\
2p(2-3p)(2-6p+5p^2),  &  \hbox{if $\lambda=2$};   \\
6p^3(3-5p),  &  \hbox{if $\lambda=3$};   \\
12p^3(2+p),  &  \hbox{if $\lambda=4$};   \\
\lambda(\lambda-1)^3p^4,  &  \hbox{if $\lambda\ge5$}.
\end{array}
\right.   \label{eq:FC-K13}   \\[5pt]
\FC_{\lambda}(\,\pathThree\,;p)
& = &
\left\{
\begin{array}{ll}
0,  &  \hbox{if $\lambda=1$};   \\
2p(2-3p)(2-6p+5p^2),  &  \hbox{if $\lambda=2$};   \\
24p^3(1-2p),  &  \hbox{if $\lambda=3$};   \\
\lambda(\lambda-1)^3p^4,  &  \hbox{if $\lambda\ge4$}.
\end{array}
\right.   \label{eq:FC-P3}   \\[5pt]
\FC_{\lambda}(\raisebox{-0.1cm}{\threePan};p)
& = &
\left\{
\begin{array}{ll}
0,  &  \hbox{if $\lambda\le2$};   \\
12p^2(1-3p+p^2),  &  \hbox{if $\lambda=3$};   \\
24p^3(1-p),  &  \hbox{if $\lambda=4$};   \\
\lambda(\lambda-1)^2(\lambda-2)p^4,  &  \hbox{if $\lambda\ge5$}.
\end{array}
\right.   \label{eq:FC-3Pan}   \\[5pt]
\FC_{\lambda}(C_4;p)
& = &
\left\{
\begin{array}{ll}
0,  &  \hbox{if $\lambda=1$};   \\
2p(2-3p)(2-6p+5p^2),  &  \hbox{if $\lambda=2$};   \\
6p^2(2-8p+9p^2),  &  \hbox{if $\lambda=3$};   \\
\lambda(\lambda-1)(\lambda^2-3\lambda+3)p^4,  &  \hbox{if $\lambda\ge4$}.
\end{array}
\right.   \label{eq:FC-C4}   \\[5pt]
\FC_{\lambda}(\raisebox{-0.1cm}{\WelshGraph};p)
& = &
\left\{
\begin{array}{ll}
0,  &  \hbox{if $\lambda\le2$};   \\
6p^2(5-26p+34p^2),  &  \hbox{if $\lambda=3$};   \\
48p^3(1-3p),  &  \hbox{if $\lambda=4$};   \\
\lambda(\lambda-1)(\lambda-2)^2p^4,  &  \hbox{if $\lambda\ge5$}.
\end{array}
\right.   \label{eq:FC-W}
\end{eqnarray}
\eopf
\end{propn}

When $G$ is bipartite and has at least one edge, $\FC_{1}(G;p)\equiv0$ by
Proposition \ref{propn:elem-properties}(a), and if $\lambda>\Delta(G)+1$ then
$\FC_{\lambda}(G;p)=P(G;\lambda)p^n$ by Proposition \ref{propn:FC-and-chrompoly}.
We next deal with the case $\lambda=2$, when
$\FC_{2}(G;p)$ can be computed efficiently.

\begin{thm}
\label{thm:FC-bipartite}
Let $G$ be a bipartite graph with at least one edge and whose components have sizes
$n_1,n_2,\ldots,n_{k(G)}$.
\[
\FC_{2}(G;p) ~=~
\left\{
\begin{array}{cl}
0,  &  \hbox{if $G$ has an isolated vertex};   \\
2^{k(G)} \displaystyle\prod_{i=1}^{k(G)} \left( (1-p)^{n_i} - (1-2p)^{n_i} \right),  &
\hbox{otherwise}.
\end{array}
\right.
\]
\end{thm}

\pf
Let $G$ be as stated.
If $G$ contains an isolated vertex, then applying Proposition \ref{propn:elem-properties}(a) to that
one-vertex component
and then Proposition \ref{propn:elem-properties}(e) to $G$, we obtain 0 for this case.

Suppose $G$ has no isolated vertex.  So each component has chromatic number~2.

Let the components of $G$ be $G_1,\ldots,G_{k(G)}$.  Consider any component $G_i$ and its
forced colouring polynomial $\FC_{2}(G_i;p)$.

Suppose no vertex in $G_i$ gets a colour in the initial random
partial $\lambda$-assignment.  This event has probability $(1-2p)^{n_i}$.
In this situation, no vertices in $G_i$ are forced, so no $\lambda$-colouring of $G_i$ (or of $G$)
is forced.  Symbolically, $\numFC(G_i,\emptyset;2)=0$.

It is therefore necessary that at least one vertex in $G_i$ be assigned a colour by the
random partial $\lambda$-assignment.
But once one vertex (say, $v$) gets a colour, all its neighbours are
forced to get the other colour, which in turn forces all \textit{their} neighbours, and so on, eventually
forcing a 2-colouring of the entire component $G_i$.  (This means that each vertex in $G_i$ is coloured
according to the parity of the length of any walk to it from $v$, in line with the standard way
of 2-colouring connected bipartite graphs.)  It is permissible for vertices in $G_i$ other
than $v$ to be coloured by the random $\lambda$-assignment, provided that the colours they are
initially given are the colours that would be determined by the forcing process starting with $v$.

We see from this that
if $C\not=\emptyset$ then
$\numFC(G_i,C;2)=2$, with the two 2-assignments with domain $C$ that force a 2-colouring of
$G_i$ being just the restrictions to $C$ of the two 2-colourings of $G_i$.

Applying these observations to (\ref{eq:FC-sum-over-C}), we have
\begin{eqnarray*}
\FC_2(G_i;p)  & = &
\sum_{C\subseteq V} \numFC(G_i,C;2) \, p^{|C|}(1-2p)^{n_i-|C|}   \\
& = &
\sum_{C\subseteq V : C\not=\emptyset} 2 \, p^{|C|}(1-2p)^{n_i-|C|}   \\
& = &
2 \left( \left(\sum_{C\subseteq V} p^{|C|}(1-2p)^{n_i-|C|}\right) - p^{|\emptyset|}(1-2p)^{n_i-|\emptyset|} \right)   \\
& = &
2 \left( (1-p)^{n_i} - (1-2p)^{n_i} \right) .   \\
&&
\end{eqnarray*}
Now Proposition \ref{propn:elem-properties}(e)
gives $\FC_{2}(G;p)=\prod_{i=1}^{k(G)} \FC_{2}(G_i;p)$ and the result follows. 
\eopf   \\

Note that the $\lambda=2$ cases of (\ref{eq:FC-K12}), (\ref{eq:FC-K13}), (\ref{eq:FC-P3}) and (\ref{eq:FC-C4}) are special case of this Theorem and agree
with it.

If $\chi(G)<\lambda\le\Delta(G)+1$ then the expressions for $\FC_{\lambda}(G;p)$ for
bipartite $G$ are likely to be more mathematically complex.

So far, we have discussed properties of the polynomial $\FC_{\lambda}(G;p)$ as a whole,
giving expressions and inequalities for it for various situations including some families of graphs
and a few specific graphs.  We now make a few observations about some algebraic properties of
$\FC_{\lambda}(G;p)$.

\begin{propn}
If $\lambda\ge\chi(G)$ then

(a)  the degree of $\FC_{\lambda}(G;p)$ is $n$;

(b) the multiplicity of $p$ as a factor in $\FC_{\lambda}(G;p)$ is the size of the domain of the
smallest partial $\lambda$-assignment of $G$ that forces a $\lambda$-colouring of $G$.
\end{propn}

\pf
This follows from Proposition~\ref{propn:FC-sum-over-C}.

(a) Take $C=V$ and observe that if $\lambda\ge\chi(G)$ then $\numFC(G,V;\lambda)>0$.

(b) This size is the size of the smallest
$C$ such that $\numFC(G,C;\lambda)>0$.
\eopf   \\

We will need the following property of partial $\lambda$-assignments.

\begin{propn}
\label{propn:num-nonempty-parts}
If $f:V\rightarrow[\lambda]$ is not total and $\lambda$-forces a $\lambda$-colouring of $G$ then the number of nonempty colour classes induced by $f$ is either $\lambda$ or $\lambda-1$.
\end{propn}

\pf
If $f$ is not total and $\lambda$-forces a $\lambda$-colouring of $G$ then there exists a vertex $v\not\in\dom f$
such that $f$ immediately $\lambda$-forces $v$.
But this means there are $\lambda-1$ colours appearing
among the neighbours of $v$.
So the number of nonempty colour classes must be $\lambda-1$ or $\lambda$.
\eopf   \\

The \textit{colour partition induced by $f$} is just
the set of nonempty colour classes $C(i)$ without regard to their order, i.e., $\{C(i)\mid i\in[\lambda], C(i)\not=\emptyset\}$.
Note that this is a partition of $\dom f$ rather than of $V$, so we can refer to it as
a \textit{colour partition of} $\dom f$.
A colour partition $\Pi$ is said to $\lambda$\textit{-force}, or just to \textit{force}, a $\lambda$-colouring of $G$
if every partial $\lambda$-assignment $f$ that induces $\Pi$ forces a $\lambda$-colouring of~$G$.

Write $\numFC(G,C;\lambda,k)$ for the number of partial $\lambda$-assignments with domain $C$
and $k$ nonempty colour classes that force a $\lambda$-colouring of $G$.

Write $\hbox{cp}_{\lambda}(G,C)$
for the number of colour partitions of $C$ that force
a $\lambda$-colouring of $G$.  The number of these partitions which have $k$ parts (all nonempty,
by definition) is
denoted by $\hbox{cp}_{\lambda}(G,C;k)$.
A colour partition with $k$ parts may be converted to a partial $\lambda$-assignment having those
same parts as its nonempty colour classes
in $(\lambda)_k$ ways, by elementary counting, so 
\begin{equation}
\label{eq:op-up}
\numFC(G,C;\lambda,k) = (\lambda)_k \, \hbox{cp}_{\lambda}(G,C;k) .
\end{equation}

Now consider $\numFC(G,C;\lambda)$.
If $C=V(G)$ then this is just $P(G;\lambda)$, by (\ref{eq:numFC-at-V}).
If $C\subset V(G)$ then
\begin{eqnarray*}
\numFC(G,C;\lambda)
& = &
\numFC(G,C;\lambda,\lambda-1) + \numFC(G,C;\lambda,\lambda)
~~~~~ \hbox{(by Proposition \ref{propn:num-nonempty-parts})}   \\
& = &
(\lambda)_{\lambda-1} \hbox{cp}_{\lambda}(G,C;\lambda-1) +
(\lambda)_{\lambda} \hbox{cp}_{\lambda}(G,C;\lambda)
~~~~~ \hbox{(by (\ref{eq:op-up}))}   \\
& = &
\lambda! \, ( \hbox{cp}_{\lambda}(G,C;\lambda-1) + \hbox{cp}_{\lambda}(G,C;\lambda) ) .
\end{eqnarray*}

So
\begin{equation}
\label{eq:FC-with-factorial}
\hbox{FC}_{\lambda}(G;p)\,p^n ~=~
P(G;\lambda) + \lambda! \sum_{C\subset V}  ( \hbox{cp}_{\lambda}(G,C;\lambda-1) + \hbox{cp}_{\lambda}(G,C;\lambda) ) p^{|C|}(1-\lambda p)^{n-|C|} .
\end{equation}

If $G$ has at least one edge then $P(G;\lambda)$ is a multiple of $\lambda(\lambda-1)$, so in that
case $\lambda(\lambda-1)\mid P(G;\lambda)$.  If, in addition, $\lambda=3$, then we find that
$3! \mid P(G;3)$ and therefore $3! \mid \hbox{FC}_{3}(G;p)$ by (\ref{eq:FC-with-factorial}).  Similarly, if $G$ is not
bipartite then $4! \mid \hbox{FC}_{4}(G;p)$.  More generally, we have

\begin{propn}
If $\lambda\in\{\chi(G),\chi(G)+1\}$ then $\lambda! \mid \hbox{FC}_{\lambda}(G;p)$.
\end{propn}

\pf
$P(G;\lambda)$ has the factor $(\lambda)_{\chi(G)}$, which equals $\lambda!$ if
$\chi(G)\in\{\lambda-1,\lambda\}$.  The result follows from (\ref{eq:FC-with-factorial}).
\eopf

\section{Derivatives}
\label{sec:derivatives}

We give two combinatorial interpretations for values of the derivative of $\log\FC_{\lambda}(G;p)$,
the first at $p=(\lambda+1)^{-1}$ and the second at $p=\lambda^{-1}$.

In the proofs of both of these, we will use (\ref{eq:FC-sum-over-C})
which, when differentiated with respect to $p$, becomes
\begin{eqnarray}
\frac{d}{dp} \hbox{FC}_{\lambda}(G;p)
& = &
\sum_{C\subseteq V} \numFC(G,C;\lambda) \frac{d}{dp} p^{|C|}(1-\lambda p)^{n-|C|}   \nonumber   \\
& = &
\sum_{C\subseteq V} \numFC(G,C;\lambda) \, (|C|-\lambda pn) \, p^{|C|-1}(1-\lambda p)^{n-|C|-1}  .
\label{eq:derivative-of-FC}
\end{eqnarray}

We first consider an interpretation at $p=(\lambda+1)^{-1}$.

Define $\afs(G;\lambda)$ to be the \textit{average forcing size} of $G$ for $\lambda$ colours, which
is the average size of the domain of a partial $\lambda$-assignment that forces a $\lambda$-colouring
of $G$:
\begin{equation}
\label{eq:afs-defn}
\afs(G;\lambda) ~=~
\frac{\sum_{C\subseteq V}\numFC(G,C;\lambda)\,|C|}{\sum_{C\subseteq V}\numFC(G,C;\lambda)} .
\end{equation}

\begin{thm}
\[
\afs(G;\lambda) ~=~
\frac{n\lambda}{\lambda+1} +
\frac{1}{(\lambda+1)^2} \,
\left.\frac{d}{dp} \log \FC_{\lambda}(G;p)\right|_{p=(\lambda+1)^{-1}} \,.
\]
\end{thm}

\pf
We have
\begin{eqnarray*}
\lefteqn{\left.\frac{d}{dp} \log \hbox{FC}_{\lambda}(G;p)\right|_{p=(\lambda+1)^{-1}}}   \\
& = &
\left.\left( \frac{1}{\hbox{FC}_{\lambda}(G;p)} \frac{d}{dp} \hbox{FC}_{\lambda}(G;p) \right) \right|_{p=(\lambda+1)^{-1}}   \\
& = &
\frac{\sum_{C\subseteq V}\numFC(G,C;\lambda)\left(|C|-\lambda n(\lambda+1)^{-1}\right)(\lambda+1)^{-n+2}}{\sum_{C\subseteq V}\numFC(G,C;\lambda)(\lambda+1)^{-n}}
~~~~ \hbox{(using (\ref{eq:derivative-of-FC}), (\ref{eq:FC-sum-over-C}) and (\ref{eq:prob-factor-with-one-on-lambda-plus1}))}   \\
& = &
\frac{(\lambda+1)^2\left(\sum_{C\subseteq V}\numFC(G,C;\lambda) |C|-\lambda n(\lambda+1)^{-1}\sum_{C\subseteq V}\numFC(G,C;\lambda)\right)}{\sum_{C\subseteq V}\numFC(G,C;\lambda)}   \\
& = &
(\lambda+1)^2\,\frac{\sum_{C\subseteq V}\numFC(G,C;\lambda) |C|}{\sum_{C\subseteq V}\numFC(G,C;\lambda)} - \lambda(\lambda+1)n   \\
& = &
(\lambda+1)^2\,\afs(G;\lambda) - \lambda(\lambda+1)n
~~~~ \hbox{(by (\ref{eq:afs-defn}))}.
\end{eqnarray*}
The result follows.
\eopf   \\

We now consider an interpretation at $p=\lambda^{-1}$.

A vertex $v\in V(G)$ is \textit{tight} for a partial $\lambda$-colouring $f$ of $G$
if it is immediately forced by $f$ (and is therefore not in $\dom f$).
This means that it is uncoloured by $f$ and
every colour except $\Phi^*f(v)$ appears among the neighbours of $v$.
The term is closely related, but not equivalent, to the term ``colourful'' for a vertex in a graph
colouring that has all colours except its own appearing among its neighbours.
If $v$ is not tight for $f$ then it is \textit{loose} for $f$.

Define
$\tight(G;\lambda)$ to be the
average number of tight vertices in a $\lambda$-colouring of $G$.  This may be written
\[
\tight(G;\lambda) ~=~
\frac{1}{P(G;\lambda)} \sum_{g\in\scriptCol(G;\lambda)}  ( \hbox{\# $\lambda$-tight vertices of $g$} ) .
\]

\begin{thm}
\[
\tight(G;\lambda) ~=~
\frac{n}{\lambda} -
\frac{1}{\lambda^2} \,
\left.\frac{d}{dp} \log \hbox{\rm FC}_{\lambda}(G;p)\right|_{p=\lambda^{-1}} \,.
\]
\end{thm}

\pf
From (\ref{eq:derivative-of-FC}) we obtain
\begin{eqnarray*}
\frac{d}{dp} \hbox{FC}_{\lambda}(G;p)
& = &
np^{n-1}\,\numFC(G,V;\lambda)  ~+    \\
&&
\sum_{C\subset V : |C|=n-1} \numFC(G,C;\lambda) \, (|C|-\lambda pn) \, p^{|C|-1}(1-\lambda p)^{n-|C|-1}  ~+   \\
&&
\sum_{C\subset V : |C|<n-1} \numFC(G,C;\lambda) \, (|C|-\lambda pn) \, p^{|C|-1}(1-\lambda p)^{n-|C|-1}  .
\end{eqnarray*}
For the first summand, note from (\ref{eq:numFC-at-V}) that $\numFC(G,V;\lambda) = P(G;\lambda)$.
The final sum vanishes when $p=\lambda^{-1}$, because the exponent of $1-\lambda p$ in
the last factor is always positive.
So, with some simplification in the other sum ($|C|=n-1$), we have
\begin{eqnarray*}
\left. \frac{d}{dp} \hbox{FC}_{\lambda}(G;p) \right|_{p=\lambda^{-1}}
& = &
n\lambda^{-n+1}\,P(G;\lambda)  
- \lambda^{-n+2} \sum_{C\subset V : |C|=n-1} \numFC(G,C;\lambda)    \\
& = &
n\lambda^{-n+1}\,P(G;\lambda)  
- \lambda^{-n+2} \sum_{C\subset V : |C|=n-1} ~ \sum_{g\in\scriptCol(G;\lambda)}
\llbracket\, \left.g\right|_C \,\hbox{forces}\, g \,\rrbracket    \\
& = &
n\lambda^{-n+1}\,P(G;\lambda)  
- \lambda^{-n+2} \sum_{g\in\scriptCol(G;\lambda)} ~ \sum_{C\subset V : |C|=n-1} 
\llbracket\, \left.g\right|_C \,\hbox{forces}\, g \,\rrbracket    \\
& = &
n\lambda^{-n+1}\,P(G;\lambda)  
- \lambda^{-n+2} \sum_{g\in\scriptCol(G;\lambda)} 
( \hbox{\# tight vertices of $g$} ) ,
\end{eqnarray*}
so
\begin{eqnarray*}
\left.\frac{1}{P(G;\lambda)} \frac{d}{dp} \hbox{FC}_{\lambda}(G;p)\right|_{p=\lambda^{-1}}
& = &
n\lambda^{-n+1}  
- \lambda^{-n+2} \frac{1}{P(G;\lambda)} \sum_{g\in\scriptCol(G;\lambda)} 
( \hbox{\# tight vertices of $g$} )   \\
& = &
n\lambda^{-n+1}  
- \lambda^{-n+2} \,\hbox{tight}(G;\lambda) .
\end{eqnarray*}
Therefore
\[
\left. \frac{\lambda^n}{P(G;\lambda)} \frac{d}{dp} \hbox{FC}_{\lambda}(G;p)\right|_{p=\lambda^{-1}}
~=~
n\lambda
- \lambda^{2} \,\hbox{tight}(G;\lambda) .
\]
But, using Proposition \ref{propn:elem-properties}(c), the left-hand side is just
\[
\left.\frac{1}{\hbox{FC}_{\lambda}(G;\lambda^{-1})} \frac{d}{dp} \hbox{FC}_{\lambda}(G;p)\right|_{p=\lambda^{-1}}
~=~
\left.\left( \frac{1}{\hbox{FC}_{\lambda}(G;p)} \frac{d}{dp} \hbox{FC}_{\lambda}(G;p) \right) \right|_{p=\lambda^{-1}}
~=~
\left.\frac{d}{dp} \log \hbox{FC}_{\lambda}(G;p)\right|_{p=\lambda^{-1}} \,.
\]
The result follows.
\eopf

\section{Complexity}
\label{sec:complexity}

Suppose we want to compute a representation of the entire forced colouring function
$\FC(G;p,\lambda)$ for an input graph $G$.  We outlined a simple way of representing these
functions in \S\ref{sec:basic-properties}, and that entailed representing the entire chromatic
polynomial $P(G;\lambda)$ of $G$.  It is well known that computing the chromatic polynomial is
\#P-hard, since computing the number of 3-colourings of $G$ is \#P-complete \cite{linial1986}.
So computing the entire forced colouring function is \#P-hard too.

Now suppose we have fixed $\lambda\in\mathbb{N}$ and $p\in\mathbb{R}$.  We can ask, for
each input graph $G$, what is the value of $\FC(G;p,\lambda)$, i.e., of $\FC_{\lambda}(G;p)$?
This is a family of problems, one for each $\lambda$ and $p$; we will denote individual problems
in this family by $\textsc{EvalFC}(p,\lambda)$:   \\

\noindent
$\textsc{EvalFC}(p,\lambda)$   \\
\textsc{Input:}~~ Graph $G$.  \\
\textsc{Output:}~~ the value of $\FC_{\lambda}(G;p)$.   \\

Each problem in this family can be solved in exponential time by the recursive method discussed
in \cite[\S6]{farr-morgan2025}.  But can we do better?

If $\lambda=1$ then $\textsc{EvalFC}(p,\lambda)$ is trivial,
and if $\lambda=2$ then it can be done in polynomial time by Theorem~\ref{thm:FC-bipartite}
(noting that $\hbox{FC}_{2}(G;p)\equiv0$ if $G$ is not bipartite by Proposition \ref{propn:elem-properties}(a)).

Our next result gives evidence that we cannot do significantly better than an exponential-time
computation, in general.

\begin{thm}
\label{thm:FC-numberP-hard}
For every integer $\lambda\ge3$ and every $p\in(0,\lambda^{-1}]$,
the problem $\textsc{EvalFC}(p,\lambda)$ is
\#P-hard.
\end{thm}

\pf
We consider two cases: $p=\lambda^{-1}$, which we deal with in the next paragraph,
and $0<p<\lambda^{-1}$, which takes up the rest of the proof.

Firstly, suppose $p=\lambda^{-1}$.  We know that
$\hbox{FC}_{\lambda}(G;\lambda^{-1})=\lambda^{-n}P(G;\lambda)$,
by Proposition \ref{propn:elem-properties}(c),
so we can easily reduce
the computation of $P(G;\lambda)$ to the computation of $\hbox{FC}_{\lambda}(G;\lambda^{-1})$
using just one call to an oracle for the latter.  Since computation of $P(G;\lambda)$ for $\lambda\ge3$
is known to be \#P-complete \cite{linial1986},
it follows that computation of $\hbox{FC}_{\lambda}(G;\lambda^{-1})$ is \#P-hard.

Suppose now that $0<p<\lambda^{-1}$.

To prove \#P-hardness, we again use a polynomial-time Turing reduction from the problem of
computing $P(G;\lambda)$.  This uses an oracle for evaluating
$\FC_{\lambda}(H;p)$ for a suitable graph $H$ constructed from $G$.
But $H$ is no longer just $G$ itself, and the reduction is now more involved.

Let $G=(V,E)$ be any graph for which we wish to determine the value of $P(G;\lambda)$.

For any $\ell\in\mathbb{N}$, the graph $L_{\ell}(G)$ is formed from $G$
by adding $\ell$ leaves to each vertex of $G$.
We use the following names in this construction.  For each $v\in V=V(G)$, we have new
vertices $t_{v,i}$ and new edges $vt_{v,i}$ for $i\in[\ell]$.  The \textit{leafset} $T(v)$ of $v$ is
defined by $T(v):=\{t_{vi}\mid i\in[\ell]\}$ and the \textit{leafset} of $V$ is $T(V):=\bigcup_{v\in V}T(v)$.
So
\begin{eqnarray*}
V(L_{\ell}(G))  & = &  V(G) \cup \bigcup_{v\in V(G)} T(v) ,   \\
E(L_{\ell}(G))  & = &  E(G) \cup \{ vt_{v,i} \mid v\in V(G), i\in[\ell] \} .
\end{eqnarray*}
Note that $G$ is an induced subgraph of $L_{\ell}(G)$, so any $\lambda$-colouring of $L_{\ell}(G)$
gives a $\lambda$-colouring of $G$ when it is restricted to $V(G)$.

Let $\alpha$ be any constant such that
\begin{equation}
\label{eq:alpha-lwr-bd}
\alpha > \frac{\log\lambda}{\log\lambda-\log(\lambda-2)} \,.
\end{equation}
Put $\ell=\alpha n$ and construct the graph $H:=L_{\ell}(G)$.  We give $L_{\ell}(G)$ to our oracle,
which returns the value of $\hbox{FC}_{\lambda}(L_{\ell}(G);p)$.  For the rest of the proof,
we show how to use this value to obtain $P(G;\lambda)$.

Suppose $f:V(L_{\ell}(G))\rightarrow[\lambda]$ is any partial $\lambda$-assignment of $L_{\ell}(G)$.
If $f$ forces a $\lambda$-colouring of $L_{\ell}(G)$ then, since $\lambda\ge3$,
it must include all leaves in its domain:
\[
T(V)\subseteq\dom f .
\]

Let $Q(G,f)$ be the set of all vertices $v\in V(G)$ that are either in $\dom f$ or are forced by the
leaves in $T(v)$.  In the latter case, the number of colours appearing in $T(v)$ must be $\lambda-1$:
\[
Q(G,f) = \{ v\in V(G) \mid v\in\dom f \vee |f(T(v))|=\lambda-1 \} .
\]

Let $R(G,f)$ be the set of all vertices $v\in V(G)$ that are not in $\dom f$ and are not forced by $T(v)$:
\[
R(G,f) = \{ v\in V(G) \mid v\not\in\dom f \wedge |f(T(v))|\le\lambda-2 \} .
\]
Every vertex of $G$ is in at most one of these sets.  The only way a vertex in $V(G)$ can avoid being in
$Q(G,f)\cup R(G,f)$ is if there are $\lambda$ different colours appearing in $T(v)$, but in that case, $f$
is contradictory (because there is no legal colour for $v$)
and does not force a $\lambda$-colouring of $L_{\ell}(G)$.  So, if $f$ forces a $\lambda$-colouring of $G$, then $Q(G,f)$ and $R(G,f)$ partition $V(G)$.

Now let $f:V(L_{\ell}(G))\rightarrow[\lambda]$ be a random partial $\lambda$-assignment
of $L_{\ell}(G)$, chosen according to our model.
The oracle-supplied value $\FC_{\lambda}(L_{\ell}(G);p)$
gives the probability that $f$ forces a $\lambda$-colouring of $L_{\ell}(G)$,
i.e., that $\Phi^*f$ is a $\lambda$-colouring of $L_{\ell}(G)$.  When this event happens, the restriction
$\left.(\Phi^*f)\right|_{V(G)}$ must be a $\lambda$-colouring of $G$ and the sets
$Q(G,f)$ and $R(G,f)$ partition $V(G)$.  So we can partition the event
``$f$ forces a $\lambda$-colouring of $L_{\ell}(G)$'' according to the values of $Q(G,f)\in 2^{V(G)}$
and according to the $\lambda$-colouring of $G$ that is induced by $\Phi^*f$.
Taking this approach, we have
\begin{align}
\lefteqn{\hspace*{-0.5cm}\hbox{FC}_{\lambda}(L_{\ell}(G);p)}   \nonumber   \\
~= &
\Pr(\hbox{$f$ forces a $\lambda$-colouring of $L_{\ell}(G)$})   \nonumber   \\
= &
\sum_{W\subseteq V(G)} \sum_{g\in\scriptCol(G;\lambda)}
\Pr\Bigl((\hbox{$f$ forces a $\lambda$-colouring of $L_{\ell}(G)$})~\cap   \nonumber   \\[-15pt]
& ~~~~~~~~~~~~~~~~~~~~~~~~~ 
(Q(G,f)=W)\cap(R(G,f)=V(G)\setminus W)\cap(\left.\Phi^*f\right|_{V(G)}=g)\Bigr)   \nonumber   \\
= &
\sum_{W\subseteq V(G)} \sum_{g\in\scriptCol(G;\lambda)}
\Pr\Bigl((\hbox{$f$ forces a $\lambda$-colouring of $L_{\ell}(G)$})\cap(\forall v\in V(G): g(v)\not\in f(T(v)))~\cap   \nonumber   \\[-10pt]
& ~~~~~~~~~~~~~~~~~~~~~~~~~ 
(Q(G,f)=W)\cap(R(G,f)=V(G)\setminus W)\cap(\left.\Phi^*f\right|_{V(G)}=g)~\cap     \nonumber   \\
& ~~~~~~~~~~~~~~~~~~~~~~~~~ 
(T(V(G))\subseteq\dom f)\Bigr)   \nonumber   \\
&   \hbox{(since $g(v)\not\in f(T(v))$ whenever $f$ forces a $\lambda$-colouring
of $L_{\ell}(G)$ with $\Phi^*f(v)=g(v)$,}   \nonumber   \\
&   \hbox{and such an $f$ must colour all leaves)}   \nonumber   \\[5pt]
= &
\sum_{W\subseteq V(G)} \sum_{g\in\scriptCol(G;\lambda)}
\Pr\Bigl((\forall v\in V(G): g(v)\not\in f(T(v)))\cap(Q(G,f)=W)\cap(R(G,f)=V(G)\setminus W)~\cap   \nonumber   \\[-15pt]
& ~~~~~~~~~~~~~~~~~~~~~~~~~ 
(\left.\Phi^*f\right|_{V(G)}=g)\cap(T(V(G))\subseteq\dom f)\Bigr)   \nonumber   \\
&   \hbox{(since any $f$ satisfying the other conditions must force a $\lambda$-colouring of $L_{\ell}(G)$)}   \nonumber   \\[5pt]
= &
\sum_{W\subseteq V(G)} \sum_{g\in\scriptCol(G;\lambda)}
\Pr\Bigl((\forall v\in W: g(v)\not\in f(T(v)))\cap(Q(G,f)=W)~\cap   \nonumber   \\[-10pt]
& ~~~~~~~~~~~~~~~~~~~~~~~~~ 
(\left.\Phi^*f\right|_{W}=\left.g\right|_{W})\cap(T(W)\subseteq\dom f) ~\cap   \nonumber   \\[5pt]
& ~~~~~~~~~~~~~~~~~~~~~~~~~ 
(\forall v\in V(G)\setminus W: g(v)\not\in f(T(v)))\cap(R(G,f)=V(G)\setminus W)~\cap   \nonumber   \\
& ~~~~~~~~~~~~~~~~~~~~~~~~~ 
(\left.\Phi^*f\right|_{V(G)\setminus W}=\left.g\right|_{V(G)\setminus W})\cap(T(V(G)\setminus W)\subseteq\dom f)\Bigr)   
\label{eq:FC-sum-W-VminusW}   \\[5pt]
&   \hbox{(describing events according to what they mean for $W$ and $V(G)\setminus W$)} .  \nonumber   \\[5pt]
= &
\sum_{g\in\scriptCol(G;\lambda)}
\Pr\Bigl((\forall v\in V(G): g(v)\not\in f(T(v)))\cap(Q(G,f)=V(G))~\cap   \nonumber   \\[-15pt]
& ~~~~~~~~~~~~~~~~~~~~~~~~~ 
(\left.\Phi^*f\right|_{V(G)}=\left.g\right|_{V(G)})\cap(T(V(G))\subseteq\dom f) \Bigr)   \nonumber   \\[5pt]
& ~~
+~  \sum_{W\subset V(G)} \sum_{g\in\scriptCol(G;\lambda)}
\Pr\Bigl((\forall v\in W: g(v)\not\in f(T(v)))\cap(W\subseteq Q(G,f))~\cap   \nonumber   \\[-10pt]
& ~~~~~~~~~~~~~~~~~~~~~~~~~~~~~~~~~~
(\left.\Phi^*f\right|_{W}=\left.g\right|_{W})\cap(T(W)\subseteq\dom f) ~\cap   \nonumber   \\[5pt]
& ~~~~~~~~~~~~~~~~~~~~~~~~~~~~~~~~~~
(\forall v\in V(G)\setminus W: g(v)\not\in f(T(v)))\cap(V\setminus W\subseteq R(G,f))~\cap   \nonumber   \\
& ~~~~~~~~~~~~~~~~~~~~~~~~~~~~~~~~~~
(\left.\Phi^*f\right|_{V(G)\setminus W}=\left.g\right|_{V(G)\setminus W})\cap(T(V(G)\setminus W)\subseteq\dom f)\Bigr) .
\label{eq:FC-two-sums}
\end{align}
Here we have broken the sum in (\ref{eq:FC-sum-W-VminusW}) into two cases, $W=V(G)$ and $W\subset V(G)$, with simplification when $W=V$ because then $V\setminus W=\emptyset$).
Note that $Q(G,f)$ and $R(G,f)$ are constructed to be disjoint, so in the double-sum in (\ref{eq:FC-two-sums}) they must partition $V$.

We will work towards an exact expression for the first sum in (\ref{eq:FC-two-sums}),
showing it to be a multiple
of $P(G;\lambda)$, and finding an upper bound for the second sum to show that it is much smaller.
Consider the second sum, which we denote by $S_{\lambda}(L_{\ell}(G);p)$.
We drop the condition $\left.\Phi^*f\right|_{V(G)\setminus W}=\left.g\right|_{V(G)\setminus W}$ in
that sum, to get an upper bound for it.
\begin{align}
\lefteqn{S_{\lambda}(L_{\ell}(G);p)}   \nonumber   \\
\le &
\sum_{W\subset V(G)} \sum_{g\in\scriptCol(G;\lambda)}
\Pr\Bigl((\forall v\in W: g(v)\not\in f(T(v)))\cap(W\subseteq Q(G,f))~\cap   \nonumber   \\[-10pt]
& ~~~~~~~~~~~~~~~~~~~~~~~~~~~~~~~~~~
(\left.\Phi^*f\right|_{W}=\left.g\right|_{W})\cap(T(W)\subseteq\dom f) ~\cap   \nonumber   \\[5pt]
& ~~~~~~~~~~~~~~~~~~~~~~~~~~~~~~~~~~
(\forall v\in V(G)\setminus W: g(v)\not\in f(T(v)))\cap(V\setminus W\subseteq R(G,f))~\cap   \nonumber   \\
& ~~~~~~~~~~~~~~~~~~~~~~~~~~~~~~~~~~
(T(V(G)\setminus W)\subseteq\dom f)\Bigr)   
\label{eq:FC-two-sums-ineq}   \\[5pt]
= &
\sum_{W\subset V(G)} \sum_{g\in\scriptCol(G;\lambda)}
\Pr\Bigl((\forall v\in W : g(v)\not\in f(T(v)))\cap(W\subseteq Q(G,f))~\cap   \nonumber   \\[-10pt]
& ~~~~~~~~~~~~~~~~~~~~~~~~~~~~~~~~
(\left.\Phi^*f\right|_{W}=\left.g\right|_{W})\cap(T(W)\subseteq\dom f) ~\Bigr) \times   \nonumber   \\[5pt]
& ~~~~~~~~~~~~~~~~~~~~~~
\Pr\Bigl( (\forall v\in V(G)\setminus W: g(v)\not\in f(T(v)))\cap(V\setminus W\subseteq R(G,f))~\cap   \nonumber   \\
& ~~~~~~~~~~~~~~~~~~~~~~~~~~~~~~~~~~
(T(V(G)\setminus W)\subseteq\dom f)\Bigr) .
\label{eq:FC-two-sums-ineq-after-indep}
\end{align}
This last step is because the events on $W$ are now independent of the events
on $V(G)\setminus W$, due to each event now being local to the subgraph induced by a
vertex and its leafset.  This is true even of the event $\left.\Phi^*f\right|_{W}=\left.g\right|_{W}$,
which refers to forcing.  To see this, note that every $v\in W$ is either coloured by $f$
or is forced by its leafset $T(v)$, so that other vertices outside $\{v\}\cup T(v)$ (including
vertices in $V(G)\setminus W$) play no part in forcing $v$.

In order to develop our expression for $\hbox{FC}_{\lambda}(L_{\ell}(G);p)$ further,
we consider in detail the probability
\begin{eqnarray*}
\lefteqn{
\Pr\Bigl((\forall v\in W: g(v)\not\in f(T(v)))\cap(W\subseteq Q(G,f))~\cap
}      \\[-5pt]
&& ~~~~~~~~~~~~~~~~~~~~~~~~~ 
(\left.\Phi^*f\right|_{W}=\left.g\right|_{W})\cap(T(W)\subseteq\dom f)\Bigr) ,      \\
\end{eqnarray*}
which appears in the first sum of (\ref{eq:FC-two-sums}) with $W=V(G)$ and as the first factor inside
our upper bound (\ref{eq:FC-two-sums-ineq-after-indep}) for $S_{\lambda}(L_{\ell}(G);p)$,
and
\begin{eqnarray*}
\lefteqn{
\Pr\Bigl((\forall v\in V(G)\setminus W: g(v)\not\in f(T(v)))\cap(V\setminus W\subseteq R(G,f))~\cap
}      \\[-5pt]
&& ~~~~~~~~~~~~~~~~~~~~~~~~~ 
(T(V(G)\setminus W)\subseteq\dom f)\Bigr) ,      \\
\end{eqnarray*}
which appears as the second factor inside our upper bound (\ref{eq:FC-two-sums-ineq-after-indep}).

Put
\begin{eqnarray*}
\sigma(\lambda,\ell)
& := &
\sum_{i=0}^{\lambda-1}
(-1)^{i} {\lambda-1\choose i}(\lambda-1-i)^{\ell}   \\
\pi(\lambda,p,\ell)
& := &
p\cdot((\lambda-1)p)^{\ell} ~+~
(1-\lambda p)
  p^{\ell}\sigma(\lambda,\ell) ,   \\
\rho(\lambda,p,\ell)
& := &
(1-\lambda p) 
((\lambda-1) p)^{\ell} ~~ - ~~
(1-\lambda p) p^{\ell}\sigma(\lambda,\ell) .
\end{eqnarray*}
Observe that
\begin{equation}
\label{eq:pi-plus-rho}
\pi(\lambda,p,\ell) + \rho(\lambda,p,\ell)  =  
(1-(\lambda-1) p) ((\lambda-1) p)^{\ell} .
\end{equation}

The sum $\sigma(\lambda,\ell)$
has an extensive literature and history and may be expressed using Stirling numbers of the
second kind \cite{boyadzhiev2012,gould1978}.  We will not need to analyse it in detail,
but a simple lower bound will be useful later.  Observe that the terms of our sum
alternate in sign and that the sequence of their magnitudes
is strictly decreasing if $\ell$ is large enough.  In fact, $\ell>(\lambda-1)\log(\lambda-1)$ is sufficient
for this, and this is satisfied by our $\ell=\alpha n$ for large enough $n$.  Therefore we have
\begin{equation}
\label{eq:sigma-upper-bd}
\sigma(\lambda,\ell)
~\ge~
(\lambda-1)^{\ell} - (\lambda-1)(\lambda-2)^{\ell} ,
\end{equation}
which implies that
\begin{eqnarray}
\pi(\lambda,p,\ell)
& \ge &
\Bigl(p\cdot((\lambda-1)p)^{\ell}\Bigr)~+~
  \Bigl((1-\lambda p)
p^{\ell}
\bigl((\lambda-1)^{\ell} - (\lambda-1)(\lambda-2)^{\ell}\bigr)\Bigr)    \nonumber   \\
& = &
((\lambda-1)p)^{\ell}
\biggl(p~+~
(1-\lambda p)
\biggl(1- (\lambda-1)\frac{(\lambda-2)^{\ell}}{(\lambda-1)^{\ell}}\,\biggr)\biggr) .
\label{eq:stirling-type-sum-lower-bd}
\end{eqnarray}

Claim 1:

For all $W\subseteq V(G)$,
\begin{eqnarray*}
\lefteqn{
\Pr\Bigl((\forall v\in W: g(v)\not\in f(T(v)))\cap(W\subseteq Q(G,f))~\cap
}      \\[-5pt]
&& ~~~~~~~~~~~~~~~~~~~~~~~~~ 
(\left.\Phi^*f\right|_{W}=\left.g\right|_{W})\cap(T(W)\subseteq\dom f)\Bigr)      \\
& = &
\pi(\lambda,p,\ell)^{|W|}
\end{eqnarray*}

Proof of Claim 1:

Let $v\in W$.
Since $\left.\Phi^*f\right|_{W}=\left.g\right|_{W}$, the colour of $v$ under $\Phi^*f$ is $g(v)$.
If $v\in\dom f$ then $f(v)=g(v)$ and the colours in $f(T(v))$ can be any colour except $g(v)$;
they do not need to force $v$ (although they can).
If $v\not\in\dom f$ then the leaves in $T(v)$ must include all colours in $[\lambda]$
except one, since $v\in W\subseteq Q(G,f)$, and that sole excluded colour must be $g(v)$; symbolically,
$f(T(v))=[\lambda]\setminus\{g(v)\}$.

\begin{eqnarray}
\lefteqn{
\Pr\Bigl((\forall v\in W: g(v)\not\in f(T(v)))\cap(W\subseteq Q(G,f))~\cap }   \nonumber   \\[-5pt]
&& ~~~~~~~~~~~~~~~~~~~~~~~~~ 
(\left.\Phi^*f\right|_{W}=\left.g\right|_{W})\cap(T(W)\subseteq\dom f)\Bigr)   \nonumber   \\
& = &
\Pr\left(\bigcap_{v\in W}\bigl((g(v)\not\in f(T(v)))\cap(v\in Q(G,f))\cap((\Phi^*f)(v)=g(v))\cap(T(v)\subseteq\dom f)\bigr)\right)   \nonumber   \\
& = &
\Pr\Biggl(\,\bigcap_{v\in W}\biggl(\Bigl((v\in\dom f)\cap(f(v)=g(v))\cap(f(T(v))\subseteq[\lambda]\setminus\{g(v)\})\cap(T(v)\subseteq\dom f)\Bigr)~~\cup      \nonumber   \\
&&  ~~~~~~~~~~~~~~~~~
  \Bigl((v\not\in\dom f)\cap(f(T(v))=[\lambda]\setminus\{g(v)\})\cap(T(v)\subseteq\dom f)\Bigr)\biggr)\Biggr)  \nonumber   \\
& = &
\prod_{v\in W}\Pr\biggl(\Bigl((v\in\dom f)\cap(f(v)=g(v))\cap(f(T(v))\subseteq[\lambda]\setminus\{g(v)\})\cap(T(v)\subseteq\dom f)\Bigr)~~\cup      \nonumber   \\
&&  ~~~~~~~~~~~~~~
  \Bigl((v\not\in\dom f)\cap(f(T(v))=[\lambda]\setminus\{g(v)\})\cap(T(v)\subseteq\dom f)\Bigr)\biggr) .  
  \label{eq:prod-probs-long-event}
\end{eqnarray}
Now, for any $v\in V(G)$, the probability inside the product concerns two mutually exclusive events,
one within $v\in\dom f$ and the other within $v\not\in\dom f$.  We break these two events
down using independence.  We have
\begin{eqnarray*}
\Pr\Bigl((v\in\dom f)\cap(f(v)=g(v))\Bigr)  & = &  p,   \\
\Pr\Bigl((f(T(v))\subseteq[\lambda]\setminus\{g(v)\})\cap(T(v)\subseteq\dom f)\Bigr)  & = &  ((\lambda-1)p)^{\ell},
\end{eqnarray*}
with these two events being independent.  Furthermore,
\begin{eqnarray*}
\Pr(v\not\in\dom f)  & = &  1-\lambda p,
\end{eqnarray*}
and this event is independent of the event
\[
(f(T(v))=[\lambda]\setminus\{g(v)\})\cap(T(v)\subseteq\dom f) .
\]
We now consider the probability of this latter event.

For every $v\in W$,
\begin{eqnarray}
\lefteqn{\Pr\Bigl(\bigl(T(v)\subseteq\dom f\bigr)\cap\bigl(f(T(v))=[\lambda]\setminus\{g(v)\}\bigr)\Bigr)}   \nonumber   \\
& = &
\sum_{A\subseteq [\lambda]\setminus\{g(v)\}}
(-1)^{|A|}\Pr\Bigl(\bigl(T(v)\subseteq\dom f\bigr)\cap\bigl(f(T(v))\subseteq([\lambda]\setminus\{g(v)\})\setminus A\bigr)\Bigr)   \nonumber   \\[-5pt]
&&   ~~~~~~~~~~~~~~~~~~~~~~~~~~~~~~ \hbox{(by inclusion-exclusion)}   \nonumber   \\[5pt]
& = &
\sum_{A\subseteq [\lambda]\setminus\{g(v)\}}
(-1)^{|A|} \Pr\biggl(\, \bigcap_{w\in T(v)} \Bigl( \bigl(w\in\dom f\bigr)\cap\bigl( f(w)\in([\lambda]\setminus\{g(v)\})\setminus A\bigr)\Bigr)  \,\biggr)   \nonumber   \\
& = &
\sum_{A\subseteq [\lambda]\setminus\{g(v)\}}
(-1)^{|A|} \prod_{w\in T(v)} \Pr\Bigl( \bigl(w\in\dom f\bigr)\cap\bigl( f(w)\in([\lambda]\setminus\{g(v)\})\setminus A\bigr)\Bigr)   \nonumber   \\
&&   ~~~~~~~~~~~~~~~~~~~~~~ \hbox{(by independence, over $v$, of the choices for $v$ and $T(v)$)}   \nonumber   \\[5pt]
& = &
\sum_{A\subseteq [\lambda]\setminus\{g(v)\}}
(-1)^{|A|} ((\lambda-1-|A|)p)^{\ell}   \nonumber   \\
& = &
\sum_{i=0}^{\lambda-1}
(-1)^{i} {\lambda-1\choose i}((\lambda-1-i)p)^{\ell}
~~~~~~~~~~~~~~~~~~~~~~~~
\hbox{(with $i=|A|$)}   \nonumber   \\
& = &
p^{\ell}\sigma(\lambda,\ell).   \label{eq:stirling-sum}
\end{eqnarray}
Resuming at (\ref{eq:prod-probs-long-event}), and using our subsequent observations including
the expression (\ref{eq:stirling-sum}), we have
\begin{eqnarray*}
\lefteqn{
\Pr\Bigl((\forall v\in W: g(v)\not\in f(T(v)))\cap(W\subseteq Q(G,f))~\cap }   \nonumber   \\[-5pt]
&& ~~~~~~~~~~~~~~~~~~~~~~~~~ 
(\left.\Phi^*f\right|_{V(G)}=g)\cap(T(W)\subseteq\dom f)\Bigr)   \nonumber   \\
& = &
\bigl(\,p\cdot((\lambda-1)p)^{\ell} ~+~ (1-\lambda p)
p^{\ell}\sigma(\lambda,\ell)\bigr)^{|W|}   \nonumber   \\
& = &
\pi(\lambda,p,\ell)^{|W|} .
\end{eqnarray*}

This completes the proof of Claim 1.

Claim 2:
For all $W\subseteq V(G)$,
\begin{eqnarray*}
\lefteqn{
\Pr\Bigl((\forall v\in V(G)\setminus W: g(v)\not\in f(T(v)))\cap(V\setminus W\subseteq R(G,f))~\cap
}      \\[-5pt]
&& ~~~~~~~~~~~~~~~~~~~~~~~~~ 
(T(V(G)\setminus W)\subseteq\dom f)\Bigr)      \\
& = &
\rho(\lambda,p,\ell)^{n-|W|} .
\end{eqnarray*}

Proof of Claim 2:

Let $v\in V(G)\setminus W$.
Since $v\in R(G,f)$, it is uncoloured by $f$ and its leafset $T(v)$ has $\le\lambda-2$ colours,
with these all being assigned by $f$ since $T(v)\subseteq\dom f$, and none of them
being $g(v)$ since $g(v)\not\in f(T(v))$.  The probability that $v$ is uncoloured by $f$
is $1-\lambda p$, and the probability that all of $T(v)$ is coloured by $f$ and
gets $\le\lambda-2$ colours, with none
of these being $g(v)$, is
\begin{eqnarray}
\lefteqn{
\Pr\Bigl( (g(v)\not\in f(T(v)))\cap(v\in R(G,f))\cap(T(v)\subseteq\dom f) \Bigr)
}   \nonumber   \\
& = &
\Pr\Bigl( (g(v)\not\in f(T(v)))\cap(v\not\in\dom f) \cap ( |f(T(v))|\le\lambda-2)\cap(T(v)\subseteq\dom f) \Bigr)   \nonumber   \\
& = &
\Pr(v\not\in\dom f)
\Pr\Bigl( (g(v)\not\in f(T(v)))\cap( |f(T(v))|\le\lambda-2)\cap(T(v)\subseteq\dom f) \Bigr)   \nonumber   \\
& = &
(1-\lambda p)
\Pr\Bigl( (T(v)\subseteq\dom f)\cap(f(T(v))\subset[\lambda]\setminus\{g(v)\}) \Bigr)   \nonumber   \\
& = &
(1-\lambda p)
\Bigl(\,\Pr \Bigl( ( T(v)\subseteq\dom f )\cap(f(T(v))\subseteq[\lambda]\setminus\{g(v)\}) \Bigr) -   \nonumber   \\
&& ~~~~~~~~~~~~~
\Pr\Bigl( (T(v)\subseteq\dom f)\cap(f(T(v))=[\lambda]\setminus\{g(v)\}) \Bigr) \Bigr)   \nonumber   \\
& = &
(1-\lambda p) \bigl(
((\lambda-1) p)^{\ell} -
p^{\ell}\sigma(\lambda,\ell) \bigr)
~~~~~~~~~~~~~~~~~~ \hbox{(using (\ref{eq:stirling-sum})).}
\label{eq:claim2-prob-at-vertex}
\end{eqnarray}
Therefore
\begin{eqnarray*}
\lefteqn{
\Pr\Bigl((\forall v\in V(G)\setminus W: g(v)\not\in f(T(v)))\cap(V\setminus W\subseteq R(G,f))~\cap
}      \\[-5pt]
&& ~~~~~~~~~~~~~~~~~~~~~~~~~ 
(T(V(G)\setminus W)\subseteq\dom f)\Bigr)      \\
& = &
\Pr\biggl(\,\bigcap_{v\in V(G)\setminus W}
\Bigl(
(g(v)\not\in f(T(v)))\cap(v\in R(G,f))\cap(T(v)\subseteq\dom f)\Bigr)
\biggr)      \\
& = &
\prod_{v\in V(G)\setminus W}
\Pr\Bigl(
(g(v)\not\in f(T(v)))\cap(v\in R(G,f))\cap(T(v)\subseteq\dom f)\Bigr)      \\
& = &
\Bigl( (1-\lambda p)  \bigl( 
((\lambda-1) p)^{\ell} -
p^{\ell}\sigma(\lambda,\ell) \bigr) \Bigr)^{n-|W|}
~~~~~~ \hbox{(by (\ref{eq:claim2-prob-at-vertex}))}   \\
& = &
\Bigl( (1-\lambda p) 
((\lambda-1) p)^{\ell} -
(1-\lambda p) p^{\ell}\sigma(\lambda,\ell) \Bigr)^{n-|W|}    \\
& = & \rho(\lambda,p,\ell)^{n-|W|} .
\end{eqnarray*}
This completes the proof of Claim 2.   \\

By Claim 1 with $W=V$, the first sum in (\ref{eq:FC-two-sums}) is
\begin{eqnarray}
\label{eq:FC-first-sum-simplified}
\sum_{g\in\scriptCol(G;\lambda)}
\pi(\lambda,p,\ell)^n
& = &
P(G;\lambda) \,\pi(\lambda,p,\ell)^n ,
\end{eqnarray}
so by (\ref{eq:FC-two-sums}) we have
\begin{equation}
\label{eq:FC-exprn-P-pi-S}
\hbox{FC}_{\lambda}(L_{\ell}(G);p)  =  
P(G;\lambda) \,\pi(\lambda,p,\ell)^n + S_{\lambda}(L_{\ell}(G);p),
\end{equation}
where
\begin{eqnarray}
\label{eq:FC-second-sum-simplified}
S_{\lambda}(L_{\ell}(G);p)  \le
\sum_{W\subset V(G)} \sum_{g\in\scriptCol(G;\lambda)}
\pi(\lambda,p,\ell)^{|W|}
\rho(\lambda,p,\ell)^{n-|W|} ,
\end{eqnarray}
by applying Claims 1 and 2 to~(\ref{eq:FC-two-sums-ineq-after-indep}).

From (\ref{eq:FC-second-sum-simplified}),
we find simpler upper bounds for $S_{\lambda}(L_{\ell}(G);p)$.

Claim 3:
\[
S_{\lambda}(L_{\ell}(G);p)  
\le
\lambda^n
\bigl( ( \pi(\lambda,p,\ell) + \rho(\lambda,p,\ell))^{n}
- \pi(\lambda,p,\ell)^{n}
\bigr)
\]

Proof of Claim 3:
\begin{eqnarray*}
S_{\lambda}(L_{\ell}(G);p)  
& \le &
\sum_{W\subset V(G)} \sum_{g\in\scriptCol(G;\lambda)}
\pi(\lambda,p,\ell)^{|W|}
\rho(\lambda,p,\ell)^{n-|W|}   \\
& = &
P(G;\lambda) \sum_{W\subset V(G)} 
\pi(\lambda,p,\ell)^{|W|}
\rho(\lambda,p,\ell)^{n-|W|}   \\
&&  ~~~~~~~~~~~~~~~~~~~~~~~~
\hbox{(since the summands do not depend on $g$)}   \\[3pt]
& \le &
\lambda^n \sum_{W\subset V(G)} 
\pi(\lambda,p,\ell)^{|W|}
\rho(\lambda,p,\ell)^{n-|W|}   \\
& = &
\lambda^n  \biggl(\Bigl(~
\sum_{W\subseteq V(G)} 
\pi(\lambda,p,\ell)^{|W|}
\rho(\lambda,p,\ell)^{n-|W|}
\Bigr)
- \pi(\lambda,p,\ell)^{n}
~\biggr)   \\
& = &
\lambda^n
\bigl( ( \pi(\lambda,p,\ell) + \rho(\lambda,p,\ell))^{n}
- \pi(\lambda,p,\ell)^{n}
\bigr)
\end{eqnarray*}
This completes the proof of Claim 3.

Claim 4:
\[
S_{\lambda}(L_{\ell}(G);p)  <  \pi(\lambda,p,\ell)^n .
\]

Proof of Claim 4:

By our choice of $\alpha$ in (\ref{eq:alpha-lwr-bd}),
\[
\alpha\log\frac{\lambda}{\lambda-2} > \log\lambda .
\]
It follows that, for sufficiently large $n$, we have
\[
\alpha\log\frac{\lambda}{\lambda-2} > 
\frac{\log n}{n} + \log\lambda + \frac{\log(1+\lambda^{-n}) + \log\bigl(1+(n(\lambda^{n}+1))^{-1}\bigr) + \log(\lambda-1)}{n} ,
\]
and therefore, since $\ell=\alpha n$,
\begin{eqnarray*}
\ell 
& > & \Bigl( \log \Bigl( \frac{
\lambda
}{
\lambda-2
} \Bigr) \Bigr)^{-1}
\Bigl( \log n + n \log\lambda + \log(1+\lambda^{-n}) + \log\bigl(1+(n(\lambda^{n}+1))^{-1}\bigr) + \log(\lambda-1) \Bigr)   \\
& = & \Bigl( \log \Bigl( \frac{
 \lambda
}{
\lambda-2
} \Bigr) \Bigr)^{-1}
\Bigl( \log\bigl((n\lambda^{n}(1+\lambda^{-n}))(1+(n(\lambda^{n}+1))^{-1}\bigr) + \log(\lambda-1) \Bigr) .
\end{eqnarray*}
Some routine manipulation gives a series of equivalent inequalities:
\begin{eqnarray*}
\ell \log \Bigl( \frac{
 \lambda
}{
\lambda-2
} \Bigr)
& > & \log(n(\lambda^{n}+1)+1) + \log(\lambda-1) ,   \\
\Bigl( \frac{
 \lambda
}{
\lambda-2
} \Bigr)^{\ell}
& > & (n(\lambda^{n}+1)+1) (\lambda-1) ,   \\
\frac{
 \lambda^{\ell}
}{
(\lambda-1)(\lambda-2)^{\ell}
} 
& > & n(\lambda^{n}+1)+1 .
\end{eqnarray*}
Therefore, for any $p\in[0,\lambda^{-1})$,
\begin{eqnarray*}
\frac{
(p/(1-\lambda p))(\lambda-1)^{\ell} +  \lambda^{\ell}
}{
(\lambda-1)(\lambda-2)^{\ell}
} 
& > & n(\lambda^{n}+1)+1 ,
\end{eqnarray*}
which is equivalent to
\begin{eqnarray*}
\frac{
(p/(1-\lambda p))(\lambda-1)^{\ell} +  \lambda^{\ell}
}{
(\lambda-1)^{\ell} - ((\lambda-1)^{\ell}-(\lambda-1)(\lambda-2)^{\ell})
} 
& > & n(\lambda^{n}+1)+1 .
\end{eqnarray*}
Applying (\ref{eq:sigma-upper-bd}) to the subtrahend in the denominator on the left-hand side,
we have the following inequality, and then the subsequent equivalent forms of it.
\begin{eqnarray*}
\frac{
(p/(1-\lambda p))(\lambda-1)^{\ell} +  \lambda^{\ell}
}{
(\lambda-1)^{\ell} - \sigma(\lambda,\ell)
} 
& > & n(\lambda^{n}+1)+1 ,   \\
\frac{
p\cdot(\lambda-1)^{\ell} + (1-\lambda p) \lambda^{\ell}
}{
(1-\lambda p) (\lambda-1)^{\ell} - (1-\lambda p) \sigma(\lambda,\ell)
}
& > & n(\lambda^{n}+1)+1 ,   \\
n \biggl(
\frac{
p\cdot(\lambda-1)^{\ell} + (1-\lambda p) \lambda^{\ell}
}{
(1-\lambda p) (\lambda-1)^{\ell} - (1-\lambda p) \sigma(\lambda,\ell)
} -1 \biggr)^{-1}
& < & \frac{\lambda^{-n}}{1+\lambda^{-n}} .
\end{eqnarray*}
Now use the standard inequalities $x/(1+x)\le\log(1+x)\le x$ to obtain
\[
n \log \biggl( 1 + \biggl(
\frac{
p\cdot(\lambda-1)^{\ell} + (1-\lambda p) \lambda^{\ell}
}{
(1-\lambda p) (\lambda-1)^{\ell} - (1-\lambda p) \sigma(\lambda,\ell)
} -1 \biggr)^{-1}\biggr)
~<~ \log(1+\lambda^{-n})
\]
Taking the exponential function of each side, then further routine manipulation, gives
the following, with each inequality implying the next.
\begin{eqnarray*}
\biggl( 1 + \biggl(
\frac{
p\cdot(\lambda-1)^{\ell} + (1-\lambda p) \lambda^{\ell}
}{
(1-\lambda p) (\lambda-1)^{\ell} - (1-\lambda p) \sigma(\lambda,\ell)
} -1 \biggr)^{-1}\biggr)^{n}
& < & 1+\lambda^{-n} ,    \\
\biggl( 1 + \biggl(
\frac{
p\cdot(\lambda-1)^{\ell} + (1-\lambda p) \sigma(\lambda,\ell)
}{
(1-\lambda p) (\lambda-1)^{\ell} - (1-\lambda p) \sigma(\lambda,\ell)
}\biggr)^{-1}\biggr)^{n}
& < & 1+\lambda^{-n} ,    \\
\biggl( 1 + \frac{
(1-\lambda p) ((\lambda-1) p)^{\ell} - (1-\lambda p) p^{\ell}\sigma(\lambda,\ell)
}{
p\cdot((\lambda-1)p)^{\ell} + (1-\lambda p) p^{\ell}\sigma(\lambda,\ell)
}\biggr)^{n}
& < & 1+\lambda^{-n} ,    \\
\biggl( 1 + \frac{\rho(\lambda,p,\ell)}{\pi(\lambda,p,\ell)}\biggr)^{n}
& < & 1+\lambda^{-n} ,    \\
\lambda^n 
\biggl( \biggl( 1 + \frac{\rho(\lambda,p,\ell)}{\pi(\lambda,p,\ell)}\biggr)^{n}
- 1
\biggr)
& < & 1 ,    \\
\lambda^n
\bigl( ( \pi(\lambda,p,\ell) + \rho(\lambda,p,\ell))^{n}
- \pi(\lambda,p,\ell)^{n}
\bigr)   
& < & \pi(\lambda,p,\ell)^n .
\end{eqnarray*}
By Claim 3, we therefore have
\[
S_{\lambda}(L_{\ell}(G);p)  <  \pi(\lambda,p,\ell)^n .
\]
This completes the proof of Claim 4.

Now, by (\ref{eq:FC-exprn-P-pi-S}) we have
\[
\frac{\hbox{FC}_{\lambda}(L_{\ell}(G);p)}{\pi(\lambda,p,\ell)^n}  =  
P(G;\lambda) + \frac{S_{\lambda}(L_{\ell}(G);p)}{\pi(\lambda,p,\ell)^n} \,
\]
where
\[
\frac{S_{\lambda}(L_{\ell}(G);p)}{\pi(\lambda,p,\ell)^n} < 1,
\]
by Claim 4.  So $P(G;\lambda)$ can be computed as the integer part of the quotient
of the oracle-supplied $\hbox{FC}_{\lambda}(L_{\ell}(G))$ and the quantity $\pi(\lambda,p,\ell)^n$.

The Turing reduction only requires (i) computation of $L_{\ell}(G)$ for one value of $\ell$, and this
$\ell$ is linear in $n$ so the computation of $L_{\ell}(G)$ takes polynomial time;
(ii) one oracle call, to obtain $\hbox{FC}_{\lambda}(L_{\ell}(G))$;
(iii) computation of $\pi(\lambda,p,\ell)^n$, which does not depend on anything about $G$ except
its number of vertices, and takes polynomial time since $\ell$ is linear in $n$, and (iv)
at the end, one integer division of the two rational numbers obtained from (i) and (iii).
Therefore this Turing reduction is polynomial-time computable.
\eopf   \\


\section{Future work}
\label{sec:future-work}

There is much that is yet to be discovered about the polynomials $\FC_{\lambda}(G;p)$.
\begin{enumerate}[(i)]
\item  Are there any other values of $p$, outside the interval $[0,\lambda^{-1}]$, where the
value of $\FC_{\lambda}(G;p)$ has a natural combinatorial interpretation?  Candidates
to consider include $p=-1$, $p=-\lambda^{-1}$, $p=(\lambda-1)^{-1}$, $p=\frac{1}{2}$, and $p=1$.
\item  For each integer $\lambda\ge3$,
classify the real (or complex) values of $p$ according to the complexity of
$\textsc{EvalFC}(G;p,\lambda)$.
We know this is trivial for $p=0$ and \#P-hard when $0<p\le\lambda^{-1}$ (Theorem~\ref{thm:FC-numberP-hard}).
\item  Colourings are closely related to \textit{tensions} in graphs, under which (after an arbitrary
reference direction is given to each edge of a graph) values from an abelian group are assigned
to the edges so that the sum of the values on the edges around every cycle is zero (with the values
in the sum signed according to how
the direction of travel along each edge compares with the reference
direction of the edge).  Proper colourings correspond, in a natural many-to-one way, to tensions
in which no edge gets value zero.  Can forcing, and the associated counting functions, be considered
for tensions as well as for colourings?
\item  Dually, we can ask similar questions about forcing for flows in graphs.
\item  When $\lambda=-1$, the absolute value of the chromatic polynomial gives the number of
acyclic orientations of the graph \cite{stanley73}.  Is there any link between the forced colouring function and (partial?) acyclic orientations?
\end{enumerate}

Chromatic polynomials of graphs generalise readily to matroids \cite{crapo1969,rota1964}.  But the same cannot be expected
of the forced colouring function, since graphs which have the same cycle matroid can have different
forced colouring functions.  We saw this in Proposition \ref{propn:FC-some-simple-gphs}, where the two nonisomorphic trees on four vertices have different $\FC_{\lambda}(G;p)$ when $\lambda\in\{3,4\}$.

The extendable colouring function and uniquely extendable colouring functions may be worth
further investigation too, perhaps doing for them what this paper does for the forced colouring function.
But those polynomials appear to be computationally more difficult in nature, so they may be
harder to work with.

Alongside research that is specific to $\FC_{\lambda}(G;p)$, we also hope for progress
on the more general questions about graph polynomials which were raised in \cite{farr-morgan2025},
some of which are well illustrated by these forced colouring polynomials.   \\

\noindent\textbf{Acknowledgements}

I am grateful to Nina Kam\v cev, Iain Moffatt and Steven Noble for suggesting some of the questions raised above.

\end{document}